\documentclass[12pt]{amsart}
\usepackage[utf8]{inputenc}
\usepackage{amsmath}
\usepackage{amssymb}
\usepackage{amsfonts}
\usepackage{mathtools}
\usepackage{amsthm}
\usepackage{xcolor}
\numberwithin{equation}{section}

\newcommand{\F}{\mathbb{F}}
\newcommand{\Fq}{\mathbb{F}_q}
\newcommand{\Fqstar}{\mathbb{F}_q^*}

\newcommand{\Q}{\mathbb{Q}}
\newcommand{\Trk}{\operatorname{Tr}}

\newtheorem{theorem}{Theorem}

\begin{document}

\title[Distinctness of the generalized $n$-dim Kloosterman sums]
{Distinctness of the generalized $n$-dimensional Kloosterman sums}
\author{Xin Lin}
\address{Department of Mathematics, Shanghai Maritime University, Shanghai 201306, PR China.}
\email{linxin@shmtu.edu.cn}

\subjclass[2020]{11T23, 11L05}

\keywords{ Kloosterman sums, Exponential sums, Finite field, Algebraic degree, Distinctness}

\begin{abstract} In this paper, we study the algebraic degree and distinctness of a kind of generalized $n$-dimensional Kloosterman sums over finite field $\Fq$ by Galois theory and the Stickelberger theorem on $p$-adic Gauss sums.
\end{abstract}

\maketitle

\section{Introduction}

\allowdisplaybreaks[3]
	
Let $p$ be a prime number and let $\mathbb{F}_q$ denote the finite field with $q = p^r$ elements. 
Let $\zeta_p$ be a fixed primitive $p$-th root of unity $\mathbb{C}$. The nontrivial additive character $\psi: \mathbb{F}_q \to \mathbb{C}^\times$ is defined by $\psi(x) = \zeta_p^{\operatorname{Tr}(x)}$, where $\operatorname{Tr}$ is the trace map from $\mathbb{F}_{q}$ to the prime field $\mathbb{F}_p$.

Exponential sums over finite fields are fundamental objects in number theory, coding theory, and arithmetic geometry. Among these, the Kloosterman sum and its generalizations occupy a central position. To provide a unified framework covering these generalizations, let $\vec{a} = (a_1, \dots, a_n)$ and $\vec{b} = (b_1, \dots, b_n)$ be two tuples of positive integers, 
the generalized Kloosterman sum $Kl_{n, r}(\vec{a}, \vec{b}, \lambda)$ over $\mathbb{F}_q$ with parameter $\lambda \in \Fqstar$ are defined as:
\begin{equation} \label{eq:gen_kloos_def}
	Kl_{n, r}(\vec{a}, \vec{b}, \lambda) := \sum_{x_1, \dots, x_n \in \Fqstar} \psi\left( \sum_{i=1}^n x_i^{a_i} + \lambda \prod_{i=1}^n x_i^{-b_i} \right).
\end{equation}
When $\vec{a} = \vec{b} = \vec{1} =(1, \dots, 1)$, the sum \eqref{eq:gen_kloos_def} recovers the \textit{classical} $n$-dimensional Kloosterman sum $Kl_{n, r}(\lambda):=Kl_{n, r}(\vec{1}, \vec{1}, \lambda)$, for which Deligne \cite{De80} established the sharp Archimedean estimate $|Kl_{n, r}(\lambda)| \le (n+1) q^{n/2}$.

From the $p$-adic perspective, the study of these sums is remarkably rich. For the classical case $Kl_{n, r}(\lambda)$, Dwork (for $n=1$) \cite{Dw74} and Sperber (for general $n$) \cite{Sp80} analyzed the $p$-adic slopes and Newton polygons. For generalized case $Kl_{n, r}(\vec{1}, \vec{b}, \lambda)$, Wan \cite{Wan04} developed some decomposition theorems and derived conditions for the Newton polygon to attain its lower bound. The case $Kl_{n, r}(\vec{a}, \vec{1}, \lambda)$ was investigated by Bellovin et al \cite{BGO}, who described the Hodge polygons, though the slope sequences remained complex. More recently, utilizing high-dimensional Dwork cohomology, Wang and Yang \cite{WY} characterized the Hodge polygons for the general case and obtained simplified slope sequences when $\vec{b}=\vec{1}$, thereby complementing the results of Bellovin et al \cite{BGO}.

Beyond these analytic estimates, a global arithmetic problem of great interest is to view the exponential sum as an algebraic integer in $\mathbb{Q}(\zeta_p)$ and determine its algebraic degree over $\mathbb{Q}$. For the classical Kloosterman sum $Kl_{n, r}(\lambda)$, this problem was settled by Salié \cite{Sa32} for $n=k=1$. In the general $n$-dimensional case, Wan \cite{Wan95} proved that if  $p\nmid k$, the degree attains its maximal possible value:
\[
\deg(Kl_{n, r}(\lambda)) = \frac{p-1}{(n+1, p-1)}.
\]
However, the degree problem becomes subtle and remains partially open if $p \mid k$. More recently, Wan\cite{Wan21} proved that if $p > n + 2$ and $\lambda\in \F_p^*$, the degree of $Kl_{n, r}(\lambda)$ above holds as well. 
Motivated by these results, Yang \cite{Yang24} extended the degree problem to the generalized Kloosterman sums $Kl_{n}(\vec{1}, \vec{b}, \lambda)$. Under the assumptions that $p \nmid k$, $\lambda\in \F_p^*$ and all positive integers $b_i<\sqrt{\frac{p}{n}}$, Yang established that
\[
\deg(Kl_{n}(\vec{1}, \vec{b}, \lambda)) = \frac{p-1}{(\sum b_i + 1, p-1)}.
\]

To complement the existing global theory of generalized exponential sums, we consider the
following specific family of generalized Kloosterman sums in this paper
\begin{equation*}
	Kl_{n, r}(\vec{\lambda}):=
	Kl_{n, r}(\vec{a}, \vec{a}, \vec{\lambda}) = \sum_{x_1, \dots, x_n \in \Fqstar} \psi\left( \sum_{i=1}^n \lambda_i x_i^{a_i} + \lambda_{n+1} \prod_{i=1}^n x_i^{-a_i} \right),
\end{equation*}
where $\vec{\lambda} = (\lambda_1, \dots, \lambda_{n+1})$ and $\lambda_i\in\Fqstar$ for $1\leq i\leq n+1$. 
Note that if $(a_i, q-1)=1$ for all $1\leq i\leq n$, the generalized Kloosterman sums $Kl_{n, r}(\vec{a}, \vec{a}, \vec{\lambda})$ can be reduced to the classical form  $Kl_{n, r}(\vec{1}, \vec{1}, \lambda)$,  whose degree has been obtained in \cite{Sp80}. Let $d_i=(a_i, q-1)$.
Without loss of generality, we can assume that $\max_{1\leq i\leq n}d_i>1$.

Our first goal is to determine the algebraic degree of $Kl_{n,r}(\vec{\lambda})$ by $p$-adic method. 
We begin by obtaining an upper bound via Galois theory. 
\begin{theorem}\label{thm_ub}
	Let $\vec{\lambda}=\{\lambda_1,\ldots,\lambda_{n+1}\}$ , where 
	$\lambda_i\in\Fqstar$.
	Let $l$ be the least common multiple of $n$ positive integers $a_1,\ldots,a_n$. We have
	\begin{align*}
			\deg Kl_{n, r}(\vec{\lambda})
			\ \bigg|\ 
			\frac{(p-1)}{\left(\frac{l\cdot\left( n+1, p-1 \right)}{(l,q-1)}, \frac{q-1}{(l,q-1)}\right)}.
		\end{align*} 
	In particular, if $l \mid q-1$, we have
		\begin{align*}
		\deg Kl_{n, r}(\vec{\lambda})
		\ \bigg|\ 
		\frac{(p-1)}{\left(\left( n+1, p-1 \right), \frac{q-1}{l}\right)}.
	\end{align*} 
\end{theorem}

As for the lower bound, we  derive a formula for the generalized Kloosterman sums in terms of Gauss sums using Fourier inversion. Following this, we apply a finer form of Stickelberger's theorem to determine the $p$-adic main terms, and thus differentiate the Galois conjugates.


\begin{theorem}\label{thm_lb}
	Let $\vec{\lambda}=\{\lambda_1,\ldots,\lambda_{n+1}\}$ , where 
	$\lambda_i\in\Fqstar$.
	Let $a_i$ be positive integers and $d_i$ be the great common divisor of $a_i$ and $q-1$. 
	Assume $\max_{1\leq i\leq n}d_i<p-1$ and $n+1<\frac{(p-1)r }{\max_{1\leq i\leq n}d_i}$.
	If $\Trk(\lambda_1\cdots\lambda_{n+1})\neq 0$, we have
	\[\frac{p-1}{(n+1, p-1)} \ \bigg|\  \deg(Kl_{n,r}(\vec{\lambda})).\]
\end{theorem}

As a corollary, the upper bound coincides with the lower bound under some restrictions, thus providing the exact algebraic degree as follow. 
\begin{theorem}\label{thm1}
	Let $\vec{\lambda}=\{\lambda_1,\ldots,\lambda_{n+1}\}$ , where 
	$\lambda_i\in\Fqstar$.
	Let $a_i$ be positive integers, $d_i=(a_i, q-1)$ and $l=\operatorname{lcm}[a_1,\cdots,a_n]$. 
	 Assume $ l\mid q-1$ and $(n+1, p-1)$ divides $\frac{q-1}{l}$. 	Assume $\max_{1\leq i\leq n}d_i<p-1$ and $n+1<\frac{(p-1)r }{\max_{1\leq i\leq n}d_i}$.
	 If $\Trk(\lambda_1\cdots\lambda_{n+1})\neq 0$, we have
	\[\deg(Kl_{n,r}(\vec{\lambda}))=\frac{p-1}{(n+1, p-1)}.\]
\end{theorem}

The second purpose of this paper is to study the distinctness of $Kl_{n,r}(\vec{\lambda})$. 
For the classical Kloosterman sum $Kl_{n, r}(\lambda)$, Fisher \cite{Fisher} showed that if $p>(2(n+1)^{2r}+1)^2$ and $Kl_{n, r}(\lambda_1)=Kl_{n, r}(\lambda_2)$ for $\lambda_1, \lambda_2\in\Fqstar$, then $\lambda_1$ and $\lambda_2$ are Frobenius conjugate. Direct computer calculation suggested that the bound for $p$ is too restrictive. 
A referee of \cite{Fisher} conjectured a better bound for $p$ and proved the case when $q=p$. This conjecture was later extended by Wan\cite{Wan95}.
More precisely, Wan\cite{Wan95} assumed that
$p\geq(r-1)(n+1)+2$ and $p$ does not divide a certain integer. Under these restrictions, if $Kl_{n, r}(\lambda_1)=Kl_{n, r}(\lambda_2)$ for $\lambda_1, \lambda_2\in\Fqstar$, then $\lambda_1$ and $\lambda_2$ are Frobenius conjugate.
 As an application, Wan \cite{Wan95} proved that
the single Kloosterman sum $Kl_{n, r}(\lambda)$ generates the fixed field $\Q(\zeta_p)^H$ if $\Trk(\lambda)\neq 0$, where $H$ is the subgroup of $\F_p^*$ consisting of all element $t$ such that $t^{(n+1,p-1)=1}$.
Using a similar method, Zhang \cite{Zhang25} determined the generating field of the twisted Kloosterman sum $\sum_{x \in \Fqstar} \chi(x)\psi\left(x+\lambda/x \right)$ when $\Trk(\lambda)\neq 0$.

In this paper, we determine the conditions under which the $(q-1)^{n+1}$ sums $\{Kl_{n,r}(\vec{\lambda})\}$ are distinct as $\lambda_i$ ranges over $\Fqstar$. Although the method of proof is analogous to that used for the degree problem, the result obtained is stronger, with the degree formula following as a direct consequence.

\begin{theorem}\label{thm2}
	Let $\vec{\lambda}=\{\lambda_1,\ldots,\lambda_{n+1}\}$ and $\vec{\mu}=\{\mu_1,\ldots,\mu_{n+1}\}$, where 
	$\lambda_i, \mu_i\in\Fqstar$.
	Let $a_i$ be positive integers and $d_i=(a_i, q-1)$.
	Assume $\max_{1\leq i\leq n}d_i<p-1$ and $n+1<\frac{(p-1)r }{\max_{1\leq i\leq n}d_i}$.
	The generalized Kloosterman sums $Kl_{n,r}(\vec{\lambda})\neq Kl_{n,r}(\vec{\mu})$ if $\Trk(\lambda_1\cdots\lambda_{n+1}) \neq \Trk(\mu_1\cdots\mu_{n+1})$.
\end{theorem}

This paper is organized as follow. In section \ref{sec2}, we express $Kl_{n,r}(\vec{\lambda})$ in terms of the $p$-adic Gauss sums. In Section \ref{sec3} , we give the proof of Theorem \ref{thm_ub} and \ref{thm_lb}.
Theorem \ref{thm2} is proved in Section \ref{sec4}.

\section{A $p$-adic expression of $Kl_{n,r}(\vec{\lambda})$}\label{sec2}

Let $q = p^r$. Let $\mathbb{Q}_q$ be the unramified extension of $\mathbb{Q}_p$ of degree $r$ with ring of integers $\mathbb{Z}_q$. We identify the residue field of $\mathbb{Z}_q$ with the finite field $\mathbb{F}_q$. The \textit{Teichmüller character} $\omega: \Fqstar \hookrightarrow \mathbb{Z}_q^\times$ is the unique section of the reduction map satisfying $\omega(x) \equiv x \pmod{p}$. We utilize $\omega$ to identify $\Fqstar$ with the group $\mu_{q-1}$ of $(q-1)$-th roots of unity in $\mathbb{C}_p$.
Let $\chi$ be a multiplicative character on $\Fqstar$. Since $\Fqstar$ is cyclic, any such character can be uniquely written in the form $\chi = \omega^{-j}$ for an integer $j$ with $0 \leq j < q-1$, where $j=0$ corresponds to the trivial character.
The \textit{Gauss sum}  is defined by
\begin{align*}
	 g(\chi):= \sum_{u \in \Fqstar} \chi(u) \zeta_p^{\operatorname{Tr}(u)}  
	 \quad\text{or}\quad
	  g(j):= \sum_{u \in \Fqstar} \omega^{-j}(u) \zeta_p^{\operatorname{Tr}(u)}.
\end{align*}
 Note that $g(0) = -1$ and $g(pj) = g(j)$ (indices modulo $q-1$).

Crucially, the additive character values $\zeta_p^{\operatorname{Tr}(u)}$ satisfy the following Fourier inversion formula:
\begin{align*} 
	\zeta_p^{\operatorname{Tr}(u)} = \frac{1}{q-1} \sum_{j=0}^{q-2} g(j) \, \omega^j(u), \quad \text{for all } u \in \Fqstar.
\end{align*}
This relation expresses the values of the additive character as a linear combination of the Teichmüller characters, weighted by the associated Gauss sums.

The first step is to represent the generalized Kloosterman sums via Gauss sums, which will allow us to discuss its inherent structure and properties by some known theorems and results.
By the orthogonality of characters, we have
\begin{align}\label{kl}
	Kl_{n,r}(\vec{\lambda}) 
	&= \sum_{\substack{x_1,\ldots,x_n \in \Fqstar}} \psi \left(\lambda_1 x_1^{a_1} + \cdots + \lambda_n x_n^{a_n} + \frac{\lambda_{n+1}}{x_1^{a_1} \cdots x_n^{a_n}} \right)\nonumber\\
	&=\frac{1}{q-1} \sum_{\substack{x_i,\ldots,x_{n+1} \in \Fqstar}} \psi \left( \lambda_1 x_1^{a_1} + \cdots + \lambda_n x_n^{a_n} + x_{n+1} \right) \nonumber\\
	&\times\sum_{\chi} \chi \left( \frac{x_1^{a_1} \cdots x_n^{a_n}x_{n+1}}{ \lambda_{n+1} } \right)\nonumber\\
	&=\frac{1}{q-1} \sum_{\chi} \chi^{-1} (\lambda_{n+1}) \left( \prod_{i=1}^n \sum_{ x_i \in \Fqstar} \psi (\lambda_i x_i^{a_i}) \chi (x_i^{a_i}) \right) G(\chi).
\end{align}

We shall now determine the possible values of the mixed exponential sums within the final expression of formula (\ref{kl}). 
Let $\chi=\omega^{-k}$, where $0 \leq k \leq q-2$. The Fourier inversion of Gauss sums  yields the following:
\begin{align}\label{s}
	\sum_{ x_i \in \Fqstar} \psi (\lambda_i x_i^{a_i}) \chi (x_i^{a_i})
	&=\sum_{j=0}^{q-2} \frac{g(j) \omega^{j}\left(\lambda_i\right)}{q-1} 
	\sum_{x_i \in \Fqstar} \omega^{a_i(j-k)} \left(x_j \right)\nonumber\\
	&=\sum_{j \in H_{i,k}} g(j)\omega^{j}\left(\lambda_i\right),
\end{align}
where $H_{i,k} = \left\{ j \in \mathbb{Z} \mid 0 \le j \le q-2, 
\ j \equiv k (\bmod \frac{q-1}{d_i})\right\}$ 
and $d_i = (a_i, q-1)$.
Note that the last equality of formula (\ref{s}) holds since
\[ \sum_{x_i \in \Fqstar} \omega^{a_i(j-k)}\left(x_i\right)= 
\begin{cases} 
0, & \text{if } (q-1) \nmid a_i (j-k), \\
q-1, & \text{if } (q-1) \mid a_i (j-k).
\end{cases} \]
Combining formula (\ref{kl}) and (\ref{s}), it follows that:
\begin{align}\label{eq04}
	Kl_{n,r}(\vec{\lambda}) = \frac{1}{q-1} \sum_{k=0}^{q-2} \omega^k(\lambda_{n+1}) g(k) \prod_{i=1}^n \sum_{j \in H_{i,k}} g(j)\omega^{j}(\lambda_i).
\end{align}

\section{Algebraic degree of $Kl_{n,r}(\vec{\lambda})$}\label{sec3}

We now proceed to prove the algebraic degree of the generalized Kloosterman sums $Kl_{n,r}(\vec{\lambda})$. 
Recall that the generalized Kloosterman sums concerned in this paper are defined by:
\begin{align*}
	Kl_{n,r}(\vec{\lambda}) = \sum_{\substack{x_1,\ldots,x_n \in \Fqstar}} \psi \left(\lambda_1 x_1^{a_1} + \cdots + \lambda_n x_n^{a_n} + \frac{\lambda_{n+1}}{x_1^{a_1} \cdots x_n^{a_n}} \right),
\end{align*}
where $\lambda_i\in\Fqstar$ and $a_i$ are positive integers.

First, we derive an upper bound for $Kl_{n,r}(\vec{\lambda})$ through Galois theory.
The Galois group of the cyclotomic field $\Q(\zeta_p)$ over $\Q$ is isomorphic to $\F_p^*$. For $c \in \F_p^*$, the corresponding automorphism $\sigma_c$ is uniquely determined by $\sigma_c(\zeta_p) = \zeta_p^c$ and the conjugates of $Kl_{n,r}(\vec{\lambda})$ are given by $\sigma_{c}(Kl_{n,r}(\vec{\lambda})) $.

To establish a direct connection between $Kl_{n, r}(\vec{\lambda})$ and its conjugates, we consider a subgroup of $\F_{p}^*$. That is, $c \in \F_{p}^*\cap(\F_{q}^*)^l$, where $l$ is the least common multiple of  $a_1,\ldots,a_n$. 
In this case, we can denote $c=t^l$ for $t \in \F_{q}^*$. Then we have $t^{(l(p-1),q-1)}=1$ and
\begin{align*}
	&\sigma_c(Kl_{n, r}(\vec{\lambda})) \\
	&= \sum_{\substack{x_i \in \mathbb{F}_q^*}} \psi\left(\lambda_1 (t^{l/a_1}x_1)^{a_1} + \cdots +\lambda_n (t^{l/a_n}x_n)^{a_n} + \frac{t^l\cdot t^{nl}\lambda_{n+1}}{(t^{l/a_1}x_1)^{a_1} \cdots (t^{l/a_n}x_n)^{a_n}}\right) \\
	&= \sum_{\substack{x_i \in \mathbb{F}_q^*}} \psi\left(\lambda_1 x_1^{a_1} + \cdots +\lambda_n x_n^{a_n} + \frac{t^{l(n+1)}\lambda_{n+1}}{x_1^{a_1} \cdots x_n^{a_n}}\right).
\end{align*}
If $t^{\left( l(n+1), q-1 \right) } =1$, then the conjugate $\sigma_c(Kl_{n, r}(\vec{\lambda}) ) = Kl_{n, r}(\vec{\lambda}) $. 
Let
\begin{align*}
	J=\left\{ c \in \F_{p}^*\cap(\F_{q}^*)^l\ \Big|\ c=t^l, \ \text{where}\ t\in \F_{q}^* \ \text{and}\ t^{\left( l(n+1), q-1 \right) } =1 \right\}.
\end{align*}
It is obvious that $Kl_{n, r}(\vec{\lambda}) \in \Q(\zeta_p)^J$ and thus the order of $J$ provides an upper bound for $\deg Kl_{n, r}(\vec{\lambda})$.
The number of elements $c\in J$ equals to the total number of such $t$ divided by the number of solutions to $t^l = c$. That is
\begin{align*}
	|J| &=  \frac{\left( l(p-1), l(n+1), q-1 \right)}{\left(l, q-1\right)}.
\end{align*}
Hence by Galois theory, we have
\begin{align}\label{eq_ub}
	\deg Kl_{n, r}(\vec{\lambda})\ \big|\ \left[ \Q(\zeta_p)^J : \Q \right]  = \frac{(p-1)}{\frac{\left( l\cdot (p-1, n+1), q-1 \right)}{\left(l, q-1\right)}}.
\end{align} 
This proves Theorem \ref{thm_ub}.

Now we shift our focus to finding an ideal lower bound, expecting that the lower bound will coincide with the upper bound as the $n+1$ parameters $\lambda_i$ varies in $\Fqstar$.
In pursuit of this, we shall determine how many distinct conjugates exist as $c$ varies in $\F_p^*$. First, we convert the conjugate element into a representation involving Gauss sums.

Similar to formula (\ref{eq04}), we acquire
\begin{align}\label{ckl}
	\sigma_c(Kl_{n,r}(\vec{\lambda})) = \frac{1}{q-1} \sum_{k=0}^{q-2} \omega^k(c\lambda_{n+1}) g(k) \prod_{i=1}^{n} \sum_{j \in H_{i,k}} g(j) \omega^j(c\lambda_i),
\end{align}
where $H_{i,k}$ are as defined in formula (\ref{s}). Then we have
\begin{align}\label{d}
	&\sigma_c(Kl_{n,r}(\vec{\lambda})) - Kl_{n,r}(\vec{\lambda}) \\
	=& \frac{1}{q-1} \sum_{k=0}^{q-2} \omega^k(\lambda_{n+1}) g(k) \nonumber\\
	&\times 
	\left(\omega^k(c) \prod_{i=1}^{n} \sum_{j \in H_{i,k}} g(j) \omega^j(c\lambda_i) - \prod_{i=1}^{n} \sum_{j \in H_{i,k}} g(j)\omega^j(\lambda_i)\right)\nonumber.
\end{align}

The next step is to discuss when the difference of Galois conjugates above is non-zero.
To address this problem, we utilize a finer form of the Stickelberger's theorem, which is an immediate consequence of the Gross-Koblitz formula \cite{GK79,Ko80}, to identify the $p$-adic leading term of the difference. The non-vanishing of the leading term then implies the non-vanishing of the difference. 

\begin{theorem}[Stickelberger]
	Let $1 \le k < q - 1$ and let $k$ have the base $p$ expansion $k = k_0 + k_1 p + \cdots + k_{r-1} p^{r-1}$. Let $\sigma(k) = k_0 + \cdots + k_{r-1}$ be the sum of the digits.
	Then, the Gauss sum $g(k)$ satisfies the congruence
	\begin{equation} \label{eq:finer_stickelberger_cong}
		g(k) \equiv - \frac{\pi^{\sigma(k)}}{k_0! \cdot k_1! \cdots k_{r-1}!} \left(\bmod{\pi^{\sigma(k) + p-1}}\right),
	\end{equation}
	where $\pi$ is the $(p-1)$-th root of $(-p)$ such that $\pi \equiv \Psi(1) - 1 (\bmod{\pi^2})$.
\end{theorem}

Recall that $H_{i,k} = \left\{ j \in \mathbb{Z} \mid 0 \le j \le q-2, \ j \equiv k (\bmod \frac{q-1}{d_i})\right\}$, where $d_i =(a_i, q-1)$ and $q=p^r$ with $r\geq 1$. 
In the following, we shall determine that as the integer $j$ and $k$ varies in formula (\ref{d}), when does the $p$-adic valuation of  $\sigma_c(Kl_{n,r}(\vec{\lambda})) - Kl_{n,r}(\vec{\lambda})$ reaches its minimum. 
To do this, we separate the  difference  into three parts as follow.
\begin{align}\label{d2}
	&(q-1)\left(\sigma_c(Kl_{n,r}(\vec{\lambda})) - Kl_{n,r}(\vec{\lambda})\right) \\
	=& 
	\left(\prod_{i=1}^{n} \sum_{j \in H_{i,0}} g(j)\omega^{j}(\lambda_i) - 
	\prod_{i=1}^{n} \sum_{j \in H_{i,0}} g(j) \omega^{j}(c\lambda_i)\right)\nonumber\\
	&+  \mathop{\sum_{k=1}^{q-2}}_{\sigma(k)=1} \omega^k(\lambda_{n+1}) g(k) 
	\left(\omega^k(c) \prod_{i=1}^{n} \sum_{j \in H_{i,k}} g(j) \omega^j(c\lambda_i) - \prod_{i=1}^{n} \sum_{j \in H_{i,k}} g(j)\omega^j(\lambda_i)\right)\nonumber\\
	&+  \mathop{\sum_{k=1}^{q-2}}_{\sigma(k)>1} \omega^k(\lambda_{n+1}) g(k) 
	\left(\omega^k(c) \prod_{i=1}^{n} \sum_{j \in H_{i,k}} g(j) \omega^j(c\lambda_i) - \prod_{i=1}^{n} \sum_{j \in H_{i,k}} g(j)\omega^j(\lambda_i)\right)\nonumber\\
	=& S_0+S_1+S_2.\nonumber
\end{align}

When $\sigma(k)=0$, we have $k=0$ and $g(0)=-1$. This implies that 
\begin{align}\label{eq34}
	&\prod_{i=1}^{n} \sum_{j \in H_{i,0}} g(j)\omega^{j}(\lambda_i)\\
	=&\prod_{i=1}^{n}
	\left(-1+\sum_{s=1}^{d_i-1}\omega^{\frac{s(q-1)}{d_i}}(\lambda_i) g\left(\frac{s(q-1)}{d_i}\right)  \right)\nonumber\\
	=&(-1)^n+(-1)^{n-1}\sum_{i=1}^{n}\sum_{s=1}^{d_i-1}\omega^{\frac{s(q-1)}{d_i}}(\lambda_i) g\left(\frac{s(q-1)}{d_i}\right) +B_1,\nonumber
\end{align}
where $B_1$ is a summation in which each term contains a product of at least two Gauss sums. Then we have
\begin{align}\label{eq30}
	S_0=&\prod_{i=1}^{n} \sum_{s=0}^{d_i-1} g\left(\frac{s(q-1)}{d_i}\right)\omega^{\frac{s(q-1)}{d_i}}(\lambda_i) \\
	&- 	\prod_{i=1}^{n} \sum_{s=0}^{d_i-1} g\left(\frac{s(q-1)}{d_i}\right) \omega^{\frac{s(q-1)}{d_i}}(\lambda_i)\omega^{\frac{s(q-1)}{d_i}}(c)\nonumber\\
	=&(-1)^{n-1}\sum_{i=1}^{n}\sum_{s=1}^{d_i-1}
	\left(1-\omega^{\frac{s(q-1)}{d_i}}(c)\right)
	\omega^{\frac{s(q-1)}{d_i}}(\lambda_i) g\left(\frac{s(q-1)}{d_i}\right) +B_2,\nonumber
\end{align}
where $B_2$ is also a summation in which each term contains a product of at least two Gauss sums.
Recall that $q=p^r$. As proved by Wan in \cite[p.1249]{Wan21}, one has 
\begin{align*}
	\sigma\left(\frac{s(q-1)}{d_i}\right)\geq \frac{p-1}{d_i}\cdot r
\end{align*}
for $1\leq s\leq d_i-1$ with the equality holding if and only if $r=1$ and $p\equiv 1\bmod d_i$.
Assume $\max_{1\leq i\leq n}d_i<p-1$. By Stickelberger Theorem, we can deduce that
\begin{align*}
	v_{\pi}\left(g\left(\frac{s(q-1)}{d_i}\right)\right)\geq \frac{p-1}{d_i}\cdot r>r\geq 1.
\end{align*}

It can be obtained from equation (\ref{eq30}) that as $q$, $d_i$ and $c$ vary, if there exists a specific set of choices for which $\omega^{\frac{q-1}{d_i}}(c) = 1$ identically for all $1\leq i\leq n$, then $S_0 = 0$. Such a case is possible. For example, when $(p-1) \mid \frac{(q-1)}{d_i}$. Otherwise, if there exists an $i$ with $1\leq i\leq n$ such that $\omega^{\frac{q-1}{d_i}}(c) \neq 1$, we have 
\begin{align}\label{eq08}
	v_{\pi}\left(S_0\right)\geq \frac{r(p-1)}{\max_{1\leq i\leq n} d_i}.
\end{align}

Now we consider the case when $\sigma(k)=1$. 
Assume $\max_{1\leq i\leq n}d_i<p-1$. Then we have $k=1,p,p^2,\ldots,p^{r-1} < \frac{q-1}{d_i}$. This implies that
\begin{align*}
	S_1=&\mathop{\sum_{k=1}^{q-2}}_{\sigma(k)=1} \omega^k(\lambda_{n+1}) g(k) 
	\left(\omega^k(c) \prod_{i=1}^{n} \sum_{j \in H_{i,k}} g(j) \omega^j(c\lambda_i) - \prod_{i=1}^{n} \sum_{j \in H_{i,k}} g(j)\omega^j(\lambda_i)\right)\\
	=&\sum_{t=0}^{r-1}\omega^{p^t}(\lambda_{n+1}) g(p^t) 
	\left(\omega^{p^t}(c) \prod_{i=1}^{n} \sum_{s=0}^{d_i-1} g\left(p^t+\frac{s(q-1)}{d_i}\right) \omega^{p^t+\frac{s(q-1)}{d_i}}(c\lambda_i) 
	\right.\\
	&\quad	\left. - \prod_{i=1}^{n} \sum_{s=0}^{d_i-1} g\left(p^t+\frac{s(q-1)}{d_i}\right)\omega^{p^t+\frac{s(q-1)}{d_i}}(\lambda_i)\right).
\end{align*}
Since
\begin{align}\label{eq35}
	&\prod_{i=1}^{n} \sum_{s=0}^{d_i-1} g\left(p^t+\frac{s(q-1)}{d_i}\right)\omega^{p^t+\frac{s(q-1)}{d_i}}(\lambda_i)\\
	=&\prod_{i=1}^{n}\left(\omega^{p^t}(\lambda_i)g\left(p^t\right)+ \sum_{s=1}^{d_i-1} \omega^{p^t+\frac{s(q-1)}{d_i}}(\lambda_i)g\left(p^t+\frac{s(q-1)}{d_i}\right)\right)\nonumber\\
	=&\prod_{i=1}^{n}\left(\omega^{p^t}(\lambda_i)g\left(p^t\right)\right)+B_3,\nonumber
\end{align}
where $B_3$ is a summation in which each term contains a product of at least a Gauss sum $g\left(p^t+\frac{s(q-1)}{d_i}\right)$ with $s\geq 1$.
Note that $\left(p, \frac{q-1}{d_i}\right)=1$. By Stickelberger theorem, we have $v_{\pi}\left(g\left(p^t\right)\right)=1$ and $v_{\pi}\left(g\left(p^t+\frac{s(q-1)}{d_i}\right)\right)>1$ if $s\geq 1$. Then we acquire
\begin{align}\label{eq31}
	S_1=&\sum_{t=0}^{r-1}\omega^{p^t}(\lambda_{n+1}) g(p^t) 
	\left(\left(\omega^{(n+1)p^t}(c)-1\right)\left(\prod_{i=1}^{n}\omega^{p^t}(\lambda_i)\right)g^n(p^t)+B_4	\right)\\
	=&\sum_{t=0}^{r-1}\left(\left(\prod_{i=1}^{n+1}\omega^{p^t}(\lambda_i)\right)\left(\omega^{(n+1)p^t}(c)-1\right)g^{n+1}(p^t)+B_5	\right),\nonumber
\end{align}
where $B_4$ and $B_5$ are summations in which each term contains a product of at least a Gauss sum $g\left(p^t+\frac{s(q-1)}{d_i}\right)$ with $s\geq 1$.

By Stickelberger theorem, we have 
\begin{align}\label{eq32}
	v_{\pi}\left(S_1\right)=n+1
\end{align}
 if $\omega^{(n+1)p^t}(c)\neq 1$. It is obvious that the $\pi$-adic valuation of the main term of $S_2$ are greater than that of $S_1$.

Assume $n+1<\frac{p-1}{\max_{1\leq i\leq n} d_i}\cdot r$.
Combining formula (\ref{eq30}) and (\ref{eq31}), we conclude that the $p$-adic leading term of $\sigma_c(Kl_{n,r}(\vec{\lambda})) - Kl_{n,r}(\vec{\lambda})$ occurs when $\sigma(k)=1$. Combining equation (\ref{d2}), (\ref{eq08}) and (\ref{eq32}), we obtain
\begin{align}\label{eq07}
	&(q-1)\left(\sigma_c(Kl_{n,r}(\vec{\lambda})) - Kl_{n,r}(\vec{\lambda})\right) \\
	& \equiv(-1)^{n+1}\sum_{t=0}^{r-1}\omega^{p^t}(\lambda_1\cdots\lambda_{n+1})\left(\omega^{(n+1)p^t}(c)-1\right) \pi^{n+1} \left(\bmod{\pi^{n+2}}\right)\nonumber.
\end{align}
Suppose $\Trk(\lambda_1\cdots\lambda_{n+1})\neq 0$. 
The formula above implies that $\sigma_c(Kl_{n,1}(\vec{\lambda})) \neq Kl_{n,1}(\vec{\lambda})$ if $\omega^{(n+1)p^t}(c)=\omega^{(n+1)}(c) \neq 1$, which means $c^{(n+1,p-1)}\neq 1$ for $c\in\F_p^*$.
This obtains the lower bound for $\deg Kl_{n,r}(\vec{\lambda})$ that
\[\frac{p-1}{(n+1, p-1)} \ \bigg|\  \deg(Kl_{n,r}(\vec{\lambda})).\]
This completes the proof of  Theorem \ref{thm_lb}.


\section{Distinctness of $Kl_{n,r}(\vec{\lambda})$}\label{sec4}
In this section, we consider the distinctness problem of the generalized Kloosterman sums. 
For $1\leq i\leq n+1$, as the $n+1$ non-zero parameters $\lambda_i$ varies in $\Fqstar$, we shall determine when does the $(q-1)^n$ values $\{ Kl_{n,r}(\lambda_1, \ldots, \lambda_{n+1}) \}$ are distinct.

Assume $\max_{1\leq i\leq n}d_i<p-1$ and $n+1<\frac{r(p-1)}{\max_{1\leq i\leq n}d_i}$.
Applying formula (\ref{eq04}), (\ref{eq34}) and (\ref{eq35}), we acquire
\begin{align*}
	Kl_{n,r}(\lambda_1, \ldots, \lambda_{n+1})
	\equiv \frac{(-1)^{n}}{q-1}+\frac{(-1)^{n+1}}{q-1} \sum_{t=0}^{r-1} \omega^{p^t}(\lambda_1\cdots\lambda_{n+1}) \pi^{n+1}\left(\bmod{\pi^{M}}\right).
\end{align*}
Similarly, for other parameters $\mu_1,\ldots,\mu_{n+1}\in\Fqstar$, we have
\begin{align*}
	Kl_{n,r}(\mu_1,\ldots,\mu_{n+1})
	\equiv \frac{(-1)^{n+1}}{q-1} \sum_{t=0}^{r-1} \omega^{p^t}(\mu_1\cdots\mu_{n+1}) \pi^{n+1}\left(\bmod{\pi^{M}}\right).
	\end{align*}
Then
$Kl_{n,r}(\lambda_1, \ldots, \lambda_{n+1}) \not\equiv Kl_{n,r}(\mu_1,\ldots,\mu_{n+1})$
if 
\[ \sum_{t=0}^{r-1} \omega^{p^t}(\lambda_1\cdots\lambda_{n+1}) \neq \sum_{t=0}^{r-1} \omega^{p^t}(\mu_1\cdots\mu_{n+1}). \]
That is, when the traces are different:
$ \Trk(\lambda_1\cdots\lambda_{n+1}) \neq \Trk(\mu_1\cdots\mu_{n+1}) $.

\section*{Acknowledgements}
The authors would like to thank Prof. Daqing Wan for his valuable suggestions on the proof of this paper. This work is supported by National Natural Science Foundation of China(No. 12301010).

\end{document}